\documentclass[10pt]{article}

\usepackage{amsmath}
\usepackage{amssymb}
\usepackage{indentfirst}
\usepackage{graphics}
\usepackage{color}

\usepackage[normalem]{ulem}

\makeatletter

\@addtoreset{equation}{section}
\makeatother
\def\R{\mbox{\boldmath $R$}}

\def\<{\langle}
\def\>{\rangle}

\newtheorem{lem}{Lemma}[section]
\newtheorem{theo}{Theorem}[section]
\newtheorem{rem}{Remark}[section]
\newtheorem{pro}{Proposition}[section]

\newtheorem{cor}{Corollary}[section]
\makeatletter
   
   \@addtoreset{equation}{section}
\makeatother

\newcommand{\Fh}{{\mathrm{_2 F_1}}}

\begin{document}
\title{\bf Asymptotic profiles for a wave equation \\ with parameter dependent logarithmic damping}
\author{Ruy Coimbra Char\~ao\thanks{Corresponding author: ruy.charao@ufsc.br;  ruycharao@gmail.com}  \\{\small Department of Mathematics} \\{\small Federal University of Santa Catarina} \\ {\small 88040-270, Florianopolis, Brazil}\\and \\Marcello D'Abbicco\thanks{marcello.dabbicco@uniba.it} \\{\small Department of Mathematics} \\{\small University of Bari} \\{\small 70125 Bari, Italy} \\and\\Ryo Ikehata\thanks{ ikehatar@hiroshima-u.ac.jp} \\ {\small Department of Mathematics}\\ {\small Division of Educational Sciences}\\ {\small Graduate School of Humanities and Social Sciences} \\ {\small Hiroshima University} \\ {\small Higashi-Hiroshima 739-8524, Japan}}
\date{}
\maketitle
\begin{abstract}
We study a nonlocal wave equation with logarithmic damping which is rather weak in the low frequency zone as compared with frequently studied strong damping case.  We consider the Cauchy problem for this model in ${\bf R}^{n}$ and we study the asymptotic profile and optimal estimates of the solutions and the total energy as $t \to \infty$ in $L^{2}$-sense. In that case some results on hypergeometric functions are useful.
\end{abstract}
\section{Introduction}
\footnote[0]{Keywords and Phrases: Wave equation; Logarithmic $\theta$-damping; asymptotic profiles; optimal $L^{2}$ and energy decay.}
\footnote[0]{2020 Mathematics Subject Classification. Primary 35L05, 35B40 Secondary 35B05, 35B45, 35C20, 35S05.}
We consider a new type of wave equation with a logarithmic damping term
\begin{equation}
u_{tt} + Au + L_{\theta}u_{t} = 0,\ \ \ (t,x)\in (0,\infty)\times {\bf R}^{n},\label{eqn}
\end{equation}
with initial data
\begin{equation}
u(0,x)= 0,\ \ u_{t}(0,x)= u_{1}(x),\ \ \ x\in{\bf R}^{n} ,\label{initial}
\end{equation}
where $u_{1}$ is chosen as
\[ u_{1} \in L^{2}({\bf R}^{n}),\]
and a parameter dependent operator
\[L_{\theta}: D(L_{\theta}) \subset L^{2}({\bf R}^{n}) \to L^{2}({\bf R}^{n})\]
is defined for each $\theta > 0 $ as follows:
\[D(L_{\theta}) := \left\{f \in L^{2}({\bf R}^{n}) \,\bigm|\,\int_{{\bf R}^{n}}(\log(1+\vert\xi\vert^{2\theta}))^{2}\vert\hat{f}(\xi)\vert^{2}d\xi < +\infty\right\},\]
and for $f \in D(L_{\theta})$
\[(L_{\theta}f) (x) := {\cal F}_{\xi\to x}^{-1}\left(\log (1+\vert\xi\vert^{2\theta})\hat{f}(\xi)\right)(x).\]
The case $\theta=1$ has been introduced by Char\~ao-Ikehata \cite{Log-damping}.

Symbolically writing, one can see that
\[L_{\theta} = \log(I+A^{\theta}),\]
where the operator $Au := -\Delta u$ for $u \in H^{2}({\bf R}^{n})$. Here, we denote the Fourier transform ${\cal F}_{x\to\xi}(f)(\xi)$ of $f(x)$ by
\[{\cal F}_{x\to\xi}(f)(\xi) = \hat{f}(\xi) := \displaystyle{\int_{{\bf R}^{n}}}e^{-ix\cdot\xi}f(x)dx\]
as usual with $i := \sqrt{-1}$, and ${\cal F}_{\xi\to x}^{-1}$ expresses its inverse Fourier transform. Since the new operator $L_{\theta}$ is constructed by a nonnegative-valued multiplication one, it is nonnegative and self-adjoint in $L^{2}({\bf R}^{n})$. Then, by a similar argument to \cite[Proposition 2.1]{ITY} based on the Lumer-Phillips Theorem one can find that the problem (1.1)-(1.2) has a unique mild solution
\[u \in C([0,\infty);H^{1}({\bf R}^{n})) \cap C^{1}([0,\infty);L^{2}({\bf R}^{n}))\]
satisfying the energy inequality
\begin{equation}\label{energy}
E_{u}(t) \leq E_{u}(0),
\end{equation}
where
\[
E_{u}(t) := \frac{1}{2}\left(\Vert u_{t}(t,\cdot)\Vert_{L^{2}}^{2} + \Vert\nabla u(t,\cdot)\Vert_{L^{2}}^{2}\right).
\]

A main topic of this paper is to find an asymptotic profile of solutions in the $L^{2}$ topology as $t \to \infty$ to problem (1.1)-(1.2), and to apply it to get the optimal rate of decay of solutions in terms of the $L^{2}$ and energy norms. It should be noticed that when one studies the asymptotic profile of solutions to problem (1.1)-(1.2) under the moment condition
\[ \int_{{\bf R}^{n}} u_1(x)\,dx \neq 0, \]
it suffices to assume that the initial amplitude satisfies $u(0,x) = 0$, without loss of generality.\\

The asymptotic profile $\nu(t,x)$ as $t \to \infty$ of the solution $u(t,x)$ to the equation (1.1) with $\theta = 1$ is already known by \cite{Log-damping}, and has been represented as
\[\nu(t,x) = \left(\int_{{\bf R}^{n}}u_{1}(x)dx\right){\cal F}_{\xi \to x}^{-1}\left((1+\vert\xi\vert^{2})^{-\frac{t}{2}}\frac{\sin(\vert\xi\vert t)}{\vert\xi\vert}\right),\]
and the dissipative structure of the solution $u(t,x)$ as $t \to \infty$ is basically dominated by the factor
\[(1+\vert\xi\vert^{2})^{-\frac{t}{2}} = e^{-\frac{t}{2}\log(1+\vert\xi\vert^{2})}.\]
In this case, it should be noticed that the behavior of the factor $e^{-\frac{t}{2}\log(1+\vert\xi\vert^{2})}$ for small $\vert\xi\vert$ is similar to the Gauss kernel $e^{-\frac{t}{2}\vert\xi\vert^{2}}$ in the Fourier space because of the fact that
\[\lim_{\vert\xi\vert \to 0}\frac{\log(1+\vert\xi\vert^{2})}{\vert\xi\vert^{2}} = 1.\]
So, the recent result due to Char\~ao-Ikehata \cite{Log-damping} is included, in a sense, in the framework of \cite{I-14}, which dealt with the equation (1.1) with $L$ replaced by $A = -\Delta$ (strong damped waves).

A similar consideration remains valid if~$\theta\in(1/2,1)$, in the sense that, due to
\begin{equation}\label{eq:theta}
\lim_{\vert\xi\vert \to 0}\frac{\log(1+\vert\xi\vert^{2\theta})}{\vert\xi\vert^{2\theta}} = 1,
\end{equation}
the results obtained for the wave model with logarithmic damping~$L_\theta u_t$ are analogous the results which may be obtained for the wave with fractional damping~$A^\theta u_t$, namely,
\begin{equation}\label{1.5}
u_{tt} + Au + A^{\theta}u_{t} = 0,\ \ \ (t,x)\in (0,\infty)\times {\bf R}^{n}.
\end{equation}
However, when~$\theta>1$, the asymptotic profile for the wave model with logarithmic damping~$L_\theta u_t$ is, in general, different from the asymptotic profile for the wave equation with fractional damping~$A^\theta u_t$. For this latter a regularity-loss structure appears (see~\cite{II}), due to the different behavior at high frequencies:
\[
\lim_{\vert\xi\vert \to \infty}\frac{\log(1+\vert\xi\vert^{2\theta})}{\vert\xi\vert^{2}} = 0,\qquad \lim_{\vert\xi\vert \to \infty}\frac{\vert\xi\vert^{2\theta}}{\vert\xi\vert^{2}} = \infty\quad  (\theta > 1).
\]
This shows a crucial difference between \eqref{1.5} and \eqref{eqn} in the high frequency region, when~$\theta>1$. The regularity-loss structure that appears in the equation \eqref{1.5} with $\theta > 1$, do not appear in~\eqref{eqn} (see Theorem 1.1).

Finally, we observe that the dissipative structure of the solution $u(t,x)$ to problem \eqref{eqn}--\eqref{initial} as $t \to \infty$ is associated basically with the function 
\[ \psi(t,\xi)= e^{-\frac{t}{2}\log(1+\vert\xi\vert^{2\theta})}= (1+\vert\xi\vert^{2\theta})^{-\frac{t}{2}},\]
and to get the exact asymptotic behavior we handle with several results from the hypergeometric functions combined with the Gautschi inequality.

Now we  mention some previous related works.\\
After two pioneering papers due to Ponce \cite{Po} and Shibata \cite{S} studying the strongly damped wave equation
\begin{equation}\label{Strong}
u_{tt} + Au + Au_{t} = 0
\end{equation}
appear from the viewpoint of $L^{p}-L^{q}$ estimates, it seems that one of the main topics on the equation \eqref{Strong} has been shifted to study the asymptotic profile and optimal rate of decay of various norms of solutions. The first trial from the asymptotic profile of solutions (as $\to \infty$) point of view has been done in abstract form by \cite{ITY}, and in concrete form by \cite{I-14}. In this connection, before \cite{ITY, I-14}, in \cite{CLI} sharp energy decay estimates of the total energy are derived by a new type of energy method in the Fourier space combined with the Haraux-Komornik inequality. Although the sharpness of the results has been already discussed in the higher dimensional case such as $n \geq 3$ in \cite{I-14}, the low dimensional case ($n = 1,2$) has been completed in the paper \cite{IO-17} at last by observing a strong singularity. Quite recently, in the papers \cite{B, BV-1, BV-2} and \cite{Mi} higher order asymptotic expansions of the solutions as $t \to \infty$ to the equation \eqref{Strong} are investigated by finding optimal rates of decay and/or blowup in infinite time.

On the other hand, a complete generalization to the structurally damped wave equation ($\alpha, \nu > 0$):
\begin{equation}\label{1.6}
u_{tt} + \alpha A^{\sigma}u + \nu A^{\delta}u_{t} = 0
\end{equation}
with a parameter $\sigma \in (\delta, 2\delta)$ can be done in the papers \cite{DEP} and \cite{IT} (only for $\sigma = 1$) from the viewpoint of capturing the leading terms of the solutions as time goes to infinity. In this connection, we have to cite a paper due to Narazaki-Reissig \cite{NR} which studies $L^1-L^1$ estimates for solutions of the equation \eqref{1.6} with $\sigma = \alpha =1$ and $\delta \in (0,1)$ (see also $L^1-L^1$ estimates in~\cite{DA17} for waves with dissipative terms and in \cite{AR} for~$\sigma$-evolution equations).

We stress that the case~$\sigma\in(0,\delta]$ in~\eqref{1.6} corresponds to a very different asymptotic profile of the solution, and in general, two different diffusive profiles for the solution may appear, see D'Abbicco-Ebert~\cite{DAE14JDE, DAE14NA, DAE16, DAE17}.

\begin{theo}\label{main-theo}\, Let $n \geq 1$, and $u_{1} \in L^{2}({\bf R}^{n})\cap L^{1,1}({\bf R}^{n})$. Assume that $\theta>1/2$ if~$n\geq2$ or~$\theta>5/8$ if~$n=1$. Then, the unique solution $u(t,x)$ to problem {\rm (1.1)}-{\rm (1.2)} satisfies
\[\left\Vert u(t,\cdot) - \left(\int_{{\bf R}^{n}}u_{1}(x)dx\right){\cal F}_{\xi\to x}^{-1}\left((1+\vert\xi\vert^{2\theta})^{-\frac{t}{2}}\frac{\sin(\vert\xi\vert t)}{\vert\xi\vert}\right)\right\Vert_{L^{2}} \leq I_{0}t^{-\frac{n}{4\theta}+\frac{(6-8\theta)_+}{4\theta}}, \quad (t \gg 1),\]
where $(x)_+=\max\{x,0\}$, and
\[I_0=\|u_1\|_{1,1}+\|u_1\|_{L^2}.\]
In the case~$n=1$ and~$\theta\in(1/2,5/8)$, the unique solution $u(t,x)$ to problem {\rm (1.1)}-{\rm (1.2)} satisfies
\[\left\Vert u(t,\cdot) - \left(\int_{{\bf R}}u_{1}(x)dx\right){\cal F}_{\xi\to x}^{-1}\left((1+\vert\xi\vert^{2\theta})^{-\frac{t}{2}}\frac{\sin(\vert\xi\vert t)}{\vert\xi\vert}\right)\right\Vert_{L^{2}} \leq I_{0}t^{\frac1\theta-\frac32}, \quad (t \gg 1).\]
\end{theo}


\begin{theo}\label{main-theo2}
\, Let $n \geq 1$, and let $u_{1} \in L^{2}({\bf R}^{n})\cap L^{1,1}({\bf R}^{n})$. Assume that $\theta >1/2$. Then, the unique solution $u(t,x)$ to problem {\rm (1.1)}-{\rm (1.2)} satisfies\\

\vspace{0.1cm}
{\rm (i)}\,\,\,$n \geq 3$ $\Rightarrow$ $C_{n}\vert P_{1}\vert t^{-\frac{n-2}{4\theta}} \leq \Vert u(t,\cdot)\Vert_{L^{2}} \leq C_{n}^{-1}I_{0}t^{-\frac{n-2}{4\theta}}$ {\rm (}$t \gg 1${\rm )},\\

\vspace{0.1cm}
{\rm (ii)}\,\,$n = 2$ $\Rightarrow$ $C_{2}\vert P_{1}\vert \sqrt{\log t} \leq \Vert u(t,\cdot)\Vert_{L^{2}} \leq C_{2}^{-1}I_{0}\sqrt{\log t}$ {\rm (}$t \gg 1${\rm )},\\
\vspace{0.1cm}

{\rm (iii)}\,$n = 1$ $\Rightarrow$ $C_{1}\vert P_{1}\vert \sqrt{t} \leq \Vert u(t,\cdot)\Vert_{L^{2}} \leq C_{1}^{-1}I_{0}\sqrt{t}$ {\rm (}$t \gg 1${\rm )},\\

\noindent
where $I_{0}$ is defined in Theorem {\rm \ref{main-theo}},
\[P_1=\int_{{\bf R}^{n}} u_1(x)\,dx,\]
and $C_{n}$ {\rm (}$n \in {\bf N}${\rm )} are constants which  depend on $\theta$ and are independent from any $t$ and initial data.
\end{theo}
\begin{rem}\,{\rm The reason why we are particularly interested about the estimates of the solution itself (not the time and/or spatial derivatives of the solution) is that the solution itself sometimes includes a kind of singularity near $0$ frequency part, and this observation clearly appears by measuring the solution itself in terms of $L^{2}$-norm, and as for the time and/or spatial derivatives of the solution we may be able to treat sometimes by another well-known method. Anyway, we want to observe how a singularity appears in the solution throughout our series of papers. It is also interesting to note that although the estimate for the $ L^2 $-norm depends on
$\theta $, for low dimension $ n = 1 \;\mbox{and} \; 2 $ the blow-up explosion rate at infinity is the same as for the case $ \theta =1$ (see \cite{Log-damping}).}
\end{rem}

Unlike the $L^2$-norm of the solution, the energy norm decays for all dimension $n$. We observe that to obtain the optimal behavior in time for the $L^2$-norm of $u_t(t,\cdot)$ we use the asymptotic profile given by \eqref{leading2} (see Proposition \ref{proposition4.5}). In particular
the following result holds.
\begin{theo}\label{main-theo3}
\, Let $n \geq 1$, and let $u_{1} \in \left(L^{2}({\bf R}^{n})\cap L^{1,1}({\bf R}^{n})\right)$. Assume that $\theta >1/2$. Then, the unique solution $u(t,x)$ to problem {\rm (1.1)}-{\rm (1.2)} satisfies\\
\vspace{0.1cm}
 $$C_{n,\theta}|P_1|\,t^{-\frac{n}{4\theta}} \leq \Vert u_t(t,\cdot)\Vert_{L^2} + \Vert\nabla u(t,\cdot)\Vert_{L^2} \leq C_{n,\theta}^{-1}I_{0}t^{-\frac{n}{4\theta}}, \quad  t \gg 1,$$
where $C_{n,\theta}$ is a positive constant, and $I_{0}$ is a constant defined in Theorem {\rm \ref{main-theo}}.
\end{theo}
\begin{rem}
{\rm The optimality of the estimates obtained in Theorems~\ref{main-theo2} and~\ref{main-theo3} hints to the possibility to compute the critical exponent of global-in-time solutions to
\begin{equation*}
u_{tt} + Au + L_{\theta}u_{t} = |u|^p,\ \ \ (t,x)\in (0,\infty)\times {\bf R}^{n}.
\end{equation*}
Following as in~\cite{DR}, it is easy to show that global-in-time energy solutions exist for initial data small in~$L^1\cap L^2$, if~$p>p_c=1+(1+2\theta)/(n-1)$, in space dimension~$n=2$, for any~$\theta>1/2$, and if~$p\in(p_c,3]$, in space dimension~$n=3$, if~$\theta\in(1/2,3/2)$. It remains open to check whether some kind of nonexistence of global-in-time solutions result may hold for~$p\leq p_c$ or, otherwise, if the existence exponent may be improved by some means (in the case of $\sigma$-evolution equation with structural damping as in~\eqref{1.6}, see~\cite{DAE20} for a different kind of exponent).}
\end{rem}

This paper is organized as follows. In section 2 we prepare several important lemmas, which will be used later. In section 3 we shall derive the asymptotic profile of the solution as $t \to \infty$, and Theorem \ref{main-theo} can be proved at a stroke. Section 4 is divided into two subsections, and in subsection 4.1, we study the upper and lower bound of the time estimates to the profile of the solution found in Section 3, and in subsection 4.2, we study the optimality of the decay rate of the total energy by deriving the leading term of the time derivative of the solution to problem \eqref{eqn}. The result in subsection 4.2 seems new in the framework of this type of equations\\

{\bf Notation.} {\small Throughout this paper, $\| \cdot\|_q$ stands for the usual $L^q({\bf R}^{n})$-norm. For simplicity of notation, in particular, we use $\| \cdot\|$ instead of $\| \cdot\|_2$. Furthermore, we denote $\Vert\cdot\Vert_{H^{l}}$ as the usual $H^{l}$-norm. 
Furthermore, we define a relation $f(t) \sim g(t)$ as $t \to \infty$ by: there exist constant $C_{j} > 0$ ($j = 1,2$) such that
\[C_{1}g(t) \leq f(t) \leq C_{2}g(t)\quad (t \gg 1).\]

We also introduce the following weighted functional spaces.
\[L^{1,\gamma}({\bf R}^{n}) := \left\{f \in L^{1}({\bf R}^{n}) \; \bigm| \; \Vert f\Vert_{1,\gamma} := \int_{{\bf R}^{n}}(1+\vert x\vert^{\gamma})\vert f(x)\vert dx < +\infty\right\}.\]
Finally, we denote the surface area of the $n$-dimensional unit ball by $\omega_{n} := \displaystyle{\int_{\vert\omega\vert = 1}}d\omega$.

}

\section{Preliminaries}

\subsection{Basic integral estimates}

Taking advantage of the theory of hypergeometric functions, we are interested in studying the asymptotic behavior as $t\to\infty$ of special integrals.
\begin{lem}\label{lem:hyp}
Let~$0\leq x_1<x_2\leq\infty$, $\mu\in\R$, and~$t\in(0,\infty)$. Also, assume that~$\mu>0$ if~$x_1=0$ and that~$t>\mu$ if~$x_2=\infty$. We consider the integral
\[ I_{\mu;x_1,x_2}(t) = \int_{x_1}^{x_2} \frac{x^{\mu-1}}{(1+x)^t}\,dx. \]
Then the following asymptotic behavior holds, independently on~$x_2\in(x_1,\infty]$:
\begin{align}
\label{eq:I0}
& \lim_{t\to\infty} t^\mu \, I_{\mu;0,x_2}(t)
    = \Gamma(\mu), \qquad \forall \mu>0,\quad \text{if~$x_1=0$,}\\
\label{eq:Inot0}
& \lim_{t\to\infty} (x_1+1)^{t-1}\,t\, I_{\mu;x_1,x_2}(t)
    = x_1^{\mu-1}\,,\quad \text{if~$x_1>0$.}
\end{align}
\end{lem}
The relation of~$I_{\mu;x_1,x_2}(t)$ with the hypergeometric functions is based on the fact that (see, e.g., \cite[3.149]{GR}):
\begin{align}
\label{eq:Fh1}
I_{\mu;0,x_2}(t)
    & = \frac{x_2^\mu}{\mu}\,\Fh(t,\mu;\mu+1;-x_2) \\
\label{eq:Fh2}
I_{\mu;x_1,\infty}(t)
    & = \frac{x_1^{\mu-t}}{t-\mu}\,\Fh(t,t-\mu;t-\mu+1;-1/x_1).
\end{align}
The definition by series of the hypergeometric functions~$\Fh(a,b;c;z)$ is:
\[ \Fh(a,b;c;z) = \sum_{n = 0}^{\infty}\frac{(a)_{n}(b)_{n}}{(c)_{n}}\frac{z^{n}}{n!}\,. \]
{\it Proof.}
We first notice that, by difference, it is sufficient to prove~\eqref{eq:Inot0} with~$x_2=\infty$, due to~$(x_2+1)^{-t} = \textit{o}((x_1+1)^{-t})$ as~$t\to\infty$. Then we remark that (see, for instance, \cite{W}):
\[ I_{\mu;0,\infty}(t) = \int_0^\infty (1+\rho)^{-t}\,\rho^{\mu-1}\,d\rho = B(\mu,t-\mu) = \frac{\Gamma(\mu)\,\Gamma(t-\mu)}{\Gamma(t)}\,, \]
where~$B$ is the Beta function and~$\Gamma$ the Gamma function. As a consequence of the Gautschi inequality, it holds
\[ \lim_{t\to\infty} \frac{t^s\,\Gamma(t-s)}{\Gamma(t)} =1, \quad \forall s\in(0,1),\]
so that
\[ \lim_{t\to\infty} t^\mu\,I_{\mu;0,\infty}(t) = \Gamma(\mu)\,. \]
By difference, now~\eqref{eq:I0} follows as a consequence of~\eqref{eq:Inot0} with~$x_2=\infty$. Using~\eqref{eq:Fh2}, thanks to the formula
\[ \Fh (a,b;c;z) = (1-z)^{c-a-b}\,\Fh(c-a,c-b;c;z), \]
we see that
\[ I_{\mu;x_1,\infty}(t) = \frac{x_1^{\mu-1}\,(x_1+1)^{1-t}}{t-\mu}\,\Fh(1-\mu,1;t+1-\mu;-1/x_1),\]
and the proof follows by noticing that~$\Fh(1-\mu,1;t+1-\mu;-1/x_1)\to1$, as~$t\to\infty$.

\hfill
$\Box$

By the change of variable
\[ \int_\eta^{\eta_2} (1+r^{2\theta})^{-t}r^{p}dr = \frac1{2\theta}\,\int_{\eta^{2\theta}}^{\eta_2^{2\theta}} (1+r)^{-t}r^{\frac{p+1}{2\theta}-1}dr, \]
we obtain the following.
\begin{cor} \label{cor2.1}
Let $\theta>0$ and $p>-1$. Then
\begin{align}
\label{eq:Iptheta}
& \lim_{t\to\infty} t^{\frac{p+1}{2\theta}}\, I_{p,\theta}(t) = \frac1{2\theta}\,\Gamma\left(\frac{p+1}{2\theta}\right), \quad \forall p>-1,\\
\label{eq:Jptheta}
& \lim_{t\to\infty} (1+\eta^{2\theta})^{t-1}\,t\,J_{\eta;p,\theta}(t) = \frac{\eta^{p+1-2\theta}}{2\theta}, \quad \forall p\in\R,\,\eta>0,
\end{align}
where
\begin{align*}
I_{p,\theta}(t)
    & = \int_{0}^{\eta_2}(1+r^{2\theta})^{-t}r^{p}dr\,,\quad \text{for some~$\eta_2\in(0,\infty]$,}\\
J_{\eta;p,\theta}(t)
    & = \int_\eta^{\eta_2} (1+r^{2\theta})^{-t}r^{p}dr\,,\quad \text{for some~$\eta_2\in(\eta,\infty]$.}
\end{align*}
\end{cor}


\subsection{ Inequalities and asymptotics }

\begin{rem}\label{rem-theta}
{\rm Let $\theta > 0$. Then, it is important to note that  the  inequality
\begin{equation}\label{ineq-teta-teta}
|\xi|^{2\theta} - \log( 1+ |\xi|^{2\theta}) \geq 0
\end{equation}
holds for all $\xi \in  {\bf R}^n$. Moreover,  for each $\theta>1/2$ there exists a number $\delta_0=\delta_0(\theta)$, $0< \delta_0 < 1$, such that the inequality
\begin{equation*}
4|\xi|^2 - \log^2( 1+ |\xi|^{2\theta}) > 0
\end{equation*}
holds for $\xi \in  {\bf R}^n$ such that $0<|\xi| \leq \delta_0$.}
\end{rem}
\begin{rem}{\rm It seems difficult to determine all signs precisely of the function $r \mapsto 4r^{2}-\log^{2}(1+r^{2\theta})$. This is a difficult point for treating general $\theta > 1/2$. Estimates for high frequency part $\vert\xi\vert \geq \delta_{0}$ can be done by quite another method.}
\end{rem}

\begin{lem}\label{ab-estim}\,Let $ \theta >1/2$. Then, the real functions  $a(\xi) $ and $b(\xi) $ given by
 \begin{equation}\label{roj}
a(\xi)=\dfrac{\log(1+|\xi|^{2\theta })}{2} \quad \text{  and  }  \quad b(\xi)= \frac{1}{2}\sqrt{4|\xi|^{2}-\log^2(1+|\xi|^{2\theta}) }
\end{equation}
are well defined for $ \xi \in {\bf R}^{n}$ such that $0 \leq |\xi| \leq \delta_0$.

\end{lem}

To study an asymptotic profile of the solution to problem \eqref{eqn}--\eqref{initial} we consider a decomposition of the Fourier transformed initial data.

\begin{rem}\label{obs1}
{\rm Using  the  Fourier transform we can get a decomposition of the initial data $\hat{u}_1$  as follows
$$\hat{u}_1(\xi)=A_1(\xi)-iB_1(\xi)+P_1,\quad \xi \in {\bf R}^n, $$
where  $P_1, A_1, B_1$  are defined by
\[P_{1} = \int_{{\bf R}^n}u_{1}(x) dx, \quad A_{1} (\xi)=\int_{{\bf R}^n}u_{1}(x)\big(1-\cos(\xi x) \big)dx, \quad B_{1}(\xi) =\int_{{\bf R}^n}u_{1}(x)\sin(\xi x)dx.\]}
\end{rem}

The next lemma according to the above decomposition appears in  Ikehata \cite{I-04}.

\begin{lem}\label{lema2.5.4}

Let $\kappa \in [0,1]$. For $u_{1} \in L^{1,\kappa}({\bf R}^{n})$ and $\xi \in {\bf R}^{n}$ it holds that
$$|A_1(\xi)|\leq K|\xi|^\kappa\|u_{1}\|_{L^{1,\kappa}} \quad \text{ and } \quad |B_1(\xi)|\leq M|\xi|^\kappa\|u_{1}\|_{L^{1,\kappa}},$$
with positive constants $K$ and $M$ depending only on $n$.\\
\end{lem}



\section{Asymptotic profiles of solutions}

The associated Cauchy problem to \eqref{eqn}-\eqref{initial} in the Fourier space is given by
\begin{align}\label{uhat}
&\hat{u}_{tt}(t,\xi) + |\xi|^2\hat{u}(t,\xi) + \log(1+  |\xi|^{2\theta})\hat{u}_t(t,\xi) = 0,\\
&\hat{u}(0,\xi)=0, \quad
\hat{u}_{t}(0,\xi)=\hat{u}_1(\xi)\nonumber.
\end{align}

The characteristics roots  $\lambda_+$ and  $\lambda_-$  of the  characteristic polynomial
$$\lambda^2+  \log(1+|\xi|^{2\theta})\lambda +|\xi|^{2}=0, \quad  \xi \in {\bf R}^n$$
associated to the equation \eqref{uhat}
 are given by
\begin{align}\label{7.1}
\lambda_{\pm} = \dfrac{-\log(1+|\xi|^{2\theta}) \pm \sqrt{\log^2(1+|\xi|^{2\theta}) - 4| \xi|^{2}}}{2}.
\end{align}
For $\theta \geq 1/2$ it should be mentioned that there is a number $\delta_0 > 0$ such that (see \eqref{rem-theta})
\[\log^2(1+|\xi|^{2\theta}) -  4| \xi|^{2} <0\]
for $\xi \in {\bf R}^n$ with $|\xi|\leq \delta_0$. Therefore the characteristics roots are complex-valued and the real part is negative for  $\xi \in {\bf R}^n,\; |\xi| \leq \delta_0$. Then we can write down $\lambda_{\pm}$ in the following form
$$\lambda_\pm = - a(\xi)\pm i b(\xi),$$
where $a(\xi)$ and $b(\xi)$ are defined by \eqref{roj} in Lemma \ref{ab-estim}. In this case the solution of the equation \eqref{uhat} is given explicitly by
$$\hat{u}(t,\xi)= \dfrac{\hat{u}_1(\xi) }{b(\xi)}\sin(b(\xi)t)e^{-a(\xi)t}$$
for  $\xi \in {\bf R}^n, |\xi|\leq \delta_0$ and $t \geq 0$.

Next, in order to find a better expression for $\hat{u}(t,\xi)$ we apply the mean value theorem to get

\begin{equation}\label{sine}
\sin\left(b(\xi)t\right)=\sin(|\xi|t)+t\left(b(\xi)-|\xi|\right)\cos(\mu(\xi)t),
\end{equation}
with
$$\mu(\xi):=\eta_{1} b(\xi)+(1-\eta_{1})|\xi|$$
for some $\eta_{1} \in (0,1)$, and

\begin{equation}\label{root}
\frac{1}{\sqrt{1-g(r)}} = 1+\frac{\log^{2}(1+r^{2\theta})}{8r^{2}}\frac{1}{\sqrt{(1-\eta_{2}g(r))^{3}}}
\end{equation}
with some $\eta_{2} \in (0,1)$, where $r := \vert\xi\vert$, and
\[g(r) := \frac{\log^{2}(1+r^{2\theta})}{4r^{2}}.\]
\noindent
The identity \eqref{root} was obtained by applying the mean value theorem to the function $$G(s)=\frac{1}{\sqrt{(1-sg(r))^{3}}}, \; 0 \leq s \leq 1.$$
Then by using Remark \ref{obs1}, \eqref{sine} and \eqref{root} $\hat{u}(t,\xi)$ can be re-written as
\begin{align}\label{express1}
\hat{u}(t,\xi) &= P_{1}e^{-a(\xi)t}\frac{\sin(tr)}{r} + P_{1}\frac{\log^{2}(1+r^{2\theta})}{8r^{3}}\frac{1}{\sqrt{(1-\eta_{2}g(r))^{3}}}e^{-a(\xi)t}\sin(tr)\\
&+ \left(\frac{A_{1}(\xi) - iB_{1}(\xi)}{b(\xi)}\right)e^{-a(\xi)t}\sin(b(\xi)t) + P_{1}t e^{-a(\xi)t}\left(\frac{b(\xi)-r}{b(\xi)}\right)\cos(\mu(\xi)t).\nonumber
\end{align}
It should be remarked that \eqref{express1} holds for small frequency parameters $\xi \in {\bf R}_{\xi}^{n}$ satisfying $\vert\xi\vert \leq \delta_{0}$.

We now introduce a candidate to be a leading term as $t \to \infty$ of the solution in the following simple form:
\begin{equation}\label{leading}
P_{1}e^{-a(\xi)t}\dfrac{\sin(|\xi|t)}{|\xi|},
\end{equation}
where $a(\xi)=\displaystyle{\frac{\log(1+|\xi|^{2\theta})}{2}}$.

Our goal in this section is to get decay estimates in time to the remainder therms defined in \eqref{express1}. To proceed with that we define the next $three$ functions which imply remainders with respect to the leading term \eqref{express1}.
 \begin{itemize}
\item[$\bullet$] $K_1(t, \xi)=\Big(\displaystyle{\frac{A_1(\xi)-iB_1(\xi)}{b(\xi)}}\Big)e^{-a(\xi)t}\sin(b(\xi)t)$;
\item[$\bullet$] $K_2(t, \xi)= P_{1}e^{-a(\xi)t}\sin(rt)\displaystyle{\frac{\log^{2}(1+r^{2\theta})}{8r^{3}}}\displaystyle{\frac{1}{\sqrt{(1-\eta_{2}g(r))^{3}}}} ,  \quad r=|\xi|>0$;
\item[$\bullet$] $K_3(t, \xi)= t\,P_{1}e^{-a(\xi)t}\left(\dfrac{b(\xi)-|\xi|}{b(\xi)}\right)\cos(\mu(\xi)t)$,
 \end{itemize}
where $a(\xi)$ and $b(\xi)$ are defined in Lemma \ref{ab-estim}. Note that using these $K_{j}(t,\xi)$ ($j = 1,2,3$) the solution $\hat{u}(t,\xi)$ to problem \eqref{uhat} can be expressed as
\begin{equation}\label{expression}
\hat{u}(t,\xi) - P_{1}e^{-a(\xi)t}\frac{\sin(tr)}{r} = \sum_{j = 1}^{3}K_{j}(t,\xi).
\end{equation}

Let us check, in fact, that $\{K_{j}(t,\xi)\}$ become error terms by using previous lemmas studied in Section 2.

First we obtain decay rates for each one of these functions on the zone of low frequency $|\xi| \ll 1$.

We begin with the estimate for $K_1(t,\xi)$.
\vspace{0.2cm}

For this function we prepare the following expression for $1/b(\xi)$ based on \eqref{root}:
\begin{equation}\label{rootb}
\frac{1}{b(\xi)} = \frac{1}{r} + \frac{\log^{2}(1+r^{2\theta})}{8r^{3}}\frac{1}{\sqrt{(1-\eta_{2}g(r))^{3}}} , \quad r=|\xi| >0.
\end{equation}
Then,
\begin{align*}
K_{1}(t,\xi) &:= \frac{A_1(\xi)-iB_1(\xi)}{\vert\xi\vert}e^{-a(\xi)t}\sin(b(\xi)t) \\
&+ \big(A_1(\xi)-iB_1(\xi)\big)\frac{\log^{2}(1+r^{2\theta})}{8r^{3}}\frac{e^{-a(\xi)t}\sin(b(\xi)t)}{\sqrt{(1-\eta_{2}g(r))^{3}}} =: K_{1,1}(t,\xi) + K_{1,2}(t,\xi).
\end{align*}
\noindent
It is easy to check the following estimate based on Lemma \ref{lema2.5.4} with $k = 1$ and Corollary \ref{cor2.1} with $p=n-1$:
\begin{align}\label{K{1,1}}
\int_{\vert\xi\vert\leq \delta}\vert K_{1,1}(t,\xi)\vert^{2}d\xi
&\leq (M+K)^{2}\Vert u_{1}\Vert_{1,1}^{2}\int_{\vert\xi\vert\leq \delta}e^{-t\log(1+\vert\xi\vert^{2\theta})}d\xi  \nonumber\\
&= \omega_n (M+K)^{2}\Vert u_{1}\Vert_{1,1}^{2}\int_0^\delta(1+r^{2\theta})^{-t}r^{n-1}dr \nonumber \\
&\leq C\omega_{n}(M+K)^{2}t^{-\frac{n}{2\theta}}\Vert u_{1}\Vert_{1,1}^{2}, \quad (t \gg 1).
\end{align}

\begin{rem}\,{\rm We note that
\[\lim_{r \to +0}\frac{\log^{2}(1+r^{2\theta })}{r^{2}} = 1\]
for  $\theta = 1/2$ and
\[\lim_{r \to +0}\frac{\log^{2}(1+r^{2\theta })}{r^{2}} = +\infty\]
for $0 \leq \theta <1/2$. Due to theses limits the corresponding cases for $\theta$ are more difficult to treat.}
\end{rem}

Now, using the important fact
\[\lim_{r \to +0}\frac{\log^{2}(1+r^{2\theta })}{r^{2}} = 0\]
for $\theta > 1/2$, we see that there is a constant $\delta $, $0<\delta \leq 1$,  such that for all $0 < r \leq \delta$ it holds that
\begin{equation}\label{bound}
g(r) = \frac{\log^{2}(1+r^{2\theta})}{4r^{2}} \leq 1/2.
\end{equation}
Then, this implies
\begin{equation}\label{defforf}
\frac{1}{\sqrt{(1-\eta_{2}g(r))^{3}}} \leq 2\sqrt{2}.
\end{equation}
Thus,  from \eqref{bound}, \eqref{defforf} and Corollary \ref{cor2.1}, together with Lemma \ref{lema2.5.4} for $k = 1$, similarly to \eqref{K{1,1}} one can also derive
\begin{align}\label{K{1,2}}
\int_{\vert\xi\vert\leq \delta}\vert K_{1,2}(t,\xi)\vert^{2}d\xi& \leq 8^{-1}(M+K)^{2}\Vert u_{1}\Vert_{1,1}^{2}\int_{\vert\xi\vert\leq \delta}\left(\frac{\log^{2}(1+r^{2\theta})}{r^{2}}\right)^{2}e^{-2ta(\xi)}d\xi \nonumber\\
&\leq \frac{1}{2}(M+K)^{2}\Vert u_{1}\Vert_{1,1}^{2}\int_{\vert\xi\vert\leq \delta}e^{-2ta(\xi)}d\xi \nonumber\\
&= \frac{1}{2}(M+K)^{2}\Vert u_{1}\Vert_{1,1}^{2}\omega_{n}\int_{0}^{\delta}(1+r^{2\theta})^{-t}r^{n-1}d\xi \nonumber\\
&\leq C\omega_{n}(M+K)^{2}\Vert u_{1}\Vert_{1,1}^{2}t^{-\frac{n}{2\theta}}, \quad (t \gg 1).
\end{align}
By combining \eqref{K{1,1}} and \eqref{K{1,2}}  we have the  following estimate for $K_{1}(t,\xi)$,
\begin{equation}\label{K{1}}
\int_{\vert\xi\vert\leq\delta}\vert K_{1}(t,\xi)\vert^{2}d\xi \leq C_{1,n}\Vert u_{1}\Vert_{1,1}^{2}t^{-\frac{n}{2\theta}}, \quad (t \gg 1).
\end{equation}
Similarly to the computations to \eqref{K{1}} and using \eqref{defforf} one can also obtain the following estimate for $K_{2}(t,\xi)$
\begin{equation}\label{K{4}}
\int_{\vert\xi\vert\leq\delta}\vert K_{2}(t,\xi)\vert^{2}d\xi \leq C_{1,n}\,|P_{1}|^{2}\, t^{-\frac{n}{2\theta}}, \quad (t \gg 1),
\end{equation}
in the case of $\theta > 3/4$, due to
\begin{equation}\label{flim}
\lim_{r \to +0}\frac{\log^{2}(1+r^{2\theta})}{r^{3}} = 0.
\end{equation}
However, if $\theta\in(1/2,3/4]$, we have a different estimate. In this case, we estimate
\begin{equation}\label{eq:thetaloss}
\left(\frac{\log^{2}(1+r^{2\theta})}{r^{3}}\right)^{2} \leq r^{8\theta-6},
\end{equation}
so that
\begin{align}\label{K{2}}
\int_{\vert\xi\vert\leq \delta}\vert K_{2}(t,\xi)\vert^{2}d\xi& \leq 8^{-1}\,|P_1|^{2}\int_{\vert\xi\vert\leq \delta} r^{8\theta-6}\, e^{-2ta(\xi)}d\xi \nonumber\\
&= 8^{-1}\,|P_1|^{2}\omega_{n}\int_{0}^\delta(1+r^{2\theta})^{-t}r^{8\theta-6+n-1}dr \nonumber\\
&\leq C\omega_{n}|P_1|^{2}\,t^{-\frac{n+8\theta-6}{2\theta}}, \quad (t \gg 1),
\end{align}
where we used that~$8\theta-6+n>0$ for any~$\theta>1/2$ if~$n\geq2$, and for any~$\theta>5/8$ if~$n=1$. On the other hand, if~$n=1$ and~$\theta\in(1/2,5/8]$, we estimate
\begin{equation}\label{eq:thetalossn1}
\sin^2(tr)\left(\frac{\log^{2}(1+r^{2\theta})}{r^{3}}\right)^{2} \leq t\,r^{8\theta-5},
\end{equation}
so that
\begin{align}\label{K{2n1}}
\int_{\vert\xi\vert\leq \delta}\vert K_{2}(t,\xi)\vert^{2}d\xi& \leq 4^{-1}\,t\,|P_1|^{2}\int_{0}^\delta(1+r^{2\theta})^{-t}r^{8\theta-5}dr \nonumber\\
&\leq C\omega_{n}|P_1|^{2}\,t^{-3+\frac2{\theta}}, \quad (t \gg 1).
\end{align}

Finally, we have to deal with the case of $K_{3}(t,\xi)$. This part is crucial in this paper.

We need suitable estimates on this term because the multiplication by $t$ is included in its definition. To do that we observe that it is not difficult to see the following expression:
\[b(\xi) - r = r\left(-\frac{\frac{\log^{2}(1+r^{2\theta})}{4r^{2}}}{1+\sqrt{1-\frac{\log^{2}(1+r^{2\theta})}{4r^{2}}}}\right) = r\left(-\frac{g(r)}{1+\sqrt{1-g(r)}}\right),  \quad (r=|\xi|\ne 0),\]
where again
\[g(r) = \frac{\log^{2}(1+r^{2\theta})}{4r^{2}}.\]
This implies
\[\frac{b(\xi)-r}{b(\xi)} = g(r)\left(\frac{-1}{1-g(r)+\sqrt{1-g(r)}}\right) =: g(r)h(r).\]
We make a next identity to gain $r^{4\theta-2}$ near $r = 0$:
\[\left\vert \frac{b(\xi)-r}{b(\xi)}\right\vert = r^{4\theta-2}\vert h(r)\vert(\frac{g(r)}{r^{4\theta-2}})\quad r \ne 0.\]
Notice that
\[\lim_{r \to +0}\vert h(r)\vert = \frac{1}{2},\]
because of the fact
\[\lim_{r \to +0}g(r) = 0\]
for $\theta > 1/2$. Furthermore, one can check that
\[\lim_{r \to +0}\frac{g(r)}{r^{4\theta-2}} = \frac{1}{4},\]
for $2\theta > 1$. Therefore, from these facts one can find $C > 0$ and $\delta_{1}, \; 0< \delta_1 \leq 1 < \delta_{0}$ such that for all $r \in (0,\delta_{1})$
\begin{equation}\label{K_{3}-1}
\left\vert \frac{b(\xi)-r}{b(\xi)}\right\vert \leq Cr^{4\theta-2}.
\end{equation}
By \eqref{K_{3}-1} and the definition of $K_{3}(t,\xi)$ one can estimate $K_{3}(t,\xi)$ as follows:
\begin{equation}\label{delta}
\int_{\vert\xi\vert\leq \delta_{1}}\vert K_{3}(t,\xi)\vert^{2}d\xi \leq \vert P_{1}\vert^{2}t^{2}\int_{\vert\xi\vert\leq \delta_{1}}r^{8\theta-4}e^{-2ta(\xi)}d\xi \leq C\vert P_{1}\vert^{2}t^{-\frac{n+8\theta-4}{2\theta}},
\end{equation}
for each $\theta \geq 1$,  where one has just Corollary \ref{cor2.1} and the definition of $a(\xi)$ in Lemma  \ref{ab-estim}.

Note that in the case when $\theta > 1/2$ we see
\[\frac{n+8\theta-4}{2\theta} > \frac{n}{2\theta}>\frac{n-2}{2\theta}.\]


\vspace{0.2cm}
Now, by summarizing above discussion one can arrived at the following crucial lemma based on \eqref{expression}, \eqref{K{1}}, \eqref{K{4}}, \eqref{delta}.

\begin{pro}\label{proposition3.1}\,Let $n \geq 2$ and $\theta >1/2$, or~$n=1$ and~$\theta>5/8$. Then, there exists a small constant $\delta_{1} \in (0,1]$ such that
\[\int_{\vert\xi\vert\leq \delta_{1}} \big|\hat{u}(t,\xi) -  P_{1}e^{-a(\xi)t}\frac{\sin(tr)}{r}\big|^{2}d\xi
\leq C\big(\vert P_{1}\vert^{2}\,t^{\frac{(6-8\theta)_+}{2\theta}} + \Vert u_{1}\Vert_{1,1}^{2}\big)t^{-\frac{n}{2\theta}},\quad (t \gg 1),\]
with some generous constant $C=C_n > 0$ depending only on  the dimension $n$. If~$n=1$ and~$\theta\in(1/2,5/8)$, then
\[\int_{\vert\xi\vert\leq \delta_{1}} \big|\hat{u}(t,\xi) -  P_{1}e^{-a(\xi)t}\frac{\sin(tr)}{r}\big|^{2}d\xi
\leq C\big(\vert P_{1}\vert^{2}\, + \Vert u_{1}\Vert_{1,1}^{2}\big)t^{\frac2\theta-3},\quad (t \gg 1),\]
\end{pro}

Next, let us prepare the so-called high frequency estimates to the $L^2$-norm  of the solution $\hat{u}(\xi,t)$.  In fact, the solution decays very fast, as usual, to the case $\theta \geq 1/2$ on the high frequency region $|\xi| \geq \delta_1$.

 We note that on this region the characteristics  roots  can be  real on part of the region
 $\{  |\xi| \geq \delta_1\}$ depending on  the size of $\theta$, as for example $\theta \geq 3$.
 Thus, we apply another method to get the precise decay rate of the $L^2$-norm of solutions on the zone of high frequency, based on the following.
 \begin{lem}\label{lem:middle}
Assume that the roots~$\lambda_\pm$ of
\[ \lambda^2 + a\lambda + b =0 \]
verify~$\Re\lambda_- \leq \Re\lambda_+<0$. Then the solution to
\[ y''+ay'+b=0, \quad y(0)=0, \quad y'(0)=y_1, \]
verifies the decay estimate
\[ |y(t)|\leq e^{\lambda_+t}\, t\,|y_1|, \qquad |y'(t)|\leq e^{\lambda_+t}\, (1+t|\lambda_-|)\,|y_1|, \]
for any~$t\geq0$.
\end{lem}
\it{Proof.}
If~$\lambda_+=\lambda_-$, then the solution is $y=te^{\lambda_+t}\,y_1$, and the proof is concluded. Otherwise, the solution is
\[ y= \frac{e^{\lambda_+t}-e^{\lambda_-t}}{\lambda_+-\lambda_-}\,y_1 = e^{\lambda_+t}\, \frac{1-e^{(\lambda_--\lambda_+)t}}{\lambda_+-\lambda_-}\,y_1 = e^{\lambda_+t}\, t\,y_1\,\int_0^1 e^{\theta(\lambda_--\lambda_+)t}\,d\theta, \]
where we used the Taylor expansion
\[ e^x= 1 + \int_0^1 x\,e^{\theta x}\,d\theta. \]
As a consequence, using~$\Re(\lambda_--\lambda_+)\leq0$, we derive $|y(t)|\leq e^{\lambda_+t}\, t\,|y_1|$. Using
\[ y'= \frac{\lambda_+e^{\lambda_+t}-\lambda_-e^{\lambda_-t}}{\lambda_+-\lambda_-}\,y_1 = e^{\lambda_+t}y_1 + \lambda_- y \]
we conclude the proof.

\hfill
$\Box$

{\rm
\begin{pro}\label{proposition3.2}
Let $n \geq 1$, $\theta \geq 1/2$ and~$\Omega_h=\{ \xi\in\R^n: \ |\xi|\geq\delta_1 \}$, where $\delta_{1} > 0$ is defined in Proposition {\rm \ref{proposition3.1}}. Then, it holds that
\[\int_{\Omega_h} |\hat{u}_t(t,\xi)|^2d\xi \leq  C\,\|u_1\|^2\, e^{ -\gamma t}, \quad(t \to \infty),\]
and
\[\int_{\Omega_h} |\hat{u}(t,\xi)|^2d\xi \leq C\,\|u_1\|^2\, e^{ -\gamma t}, \quad(t \to \infty),\]
for some constant $\gamma>0$ and $C>0$.
\end{pro}

{\it Proof.}\,By Lemma~\ref{lem:middle}, we immediately obtain
\[ \int_{\Omega_h}|\hat{u}(t,\xi)|^2d\xi \leq t^2\,e^{-2\gamma t}\,\|u_1\|^2, \]
where
\[ \gamma=\min_{\Omega_h} \Re(-\lambda_+)>0. \]
The proof follows estimating~$t^2\,e^{-\gamma t}\leq C$. To obtain the estimate for the time-derivative~$\hat u_t$, we need a preliminary step. Let~$B=B(\theta)>1$ be such that~$\log(1+r^{2\theta}) \leq r$ for any~$r\geq B$. As a consequence:
\begin{equation}\label{3.22}
4r^{2}-\log^{2}(1+r^{2\theta})\geq r\,(2r+\log(1+r^{2\theta})) \geq 2r^2
\end{equation}
for any~$r\geq B$. We divide~$\Omega_h$ into two subzones. We define
\[ \Omega_{h,1} = \{ \xi\in\R^n: \ \delta_1\leq |\xi|\leq B \},\quad \Omega_{h,2} = \{ \xi\in\R^n: \ |\xi|\geq B \}. \]
For any~$\xi\in\Omega_{h,1}$, we apply Lemma~\ref{lem:middle}. Using $|\lambda_-|^2\leq C(1+|\xi|^2) \leq C(1+B^2)$, we obtain
\[ \int_{\Omega_{h,1}}|\hat{u}_t(t,\xi)|^2d\xi \leq e^{-2\gamma t}\,\big(\|u_1\|^2+C(1+B^2)\,t^2\|u_1\|^2\big)\leq C_1\,e^{ -\gamma t}\,\|u_1\|^2. \]
On the other hand, for any~$\xi\in\Omega_{h,2}$, it holds
\[ \hat u(t,\xi) = \hat u_1(\xi)\,e^{-a(\xi)t}\,\frac{\sin b(\xi)t}{b(\xi)}, \]
where~$a(\xi)$ and~$b(\xi)$ are defined in~\eqref{roj}. Hence, we may estimate
%
\begin{align*} 
\int_{\Omega_{h,2}}|\hat{u}_t(t,\xi)|^2d\xi 
&\leq \int_{\Omega_{h,2}}\left( 1+\vert\frac{a(\xi)}{b(\xi)}\vert^{2}\right)^2 e^{-2a(\xi)t}\,|\hat u_1(\xi)|^2d\xi \\
&\leq 4 \int_{\Omega_{h,2}} e^{-2a(\xi)t}\,|\hat u_1(\xi)|^2d\xi 
\leq 4e^{-2\gamma t}\,\|u_1\|^2.
\end{align*}
Here, one has just used the fact that \eqref{3.22} implies
\[\left\vert\frac{a(\xi)}{b(\xi)}\right\vert^{2} \leq \frac{\log^{2}(1+\vert\xi\vert^{2\theta})}{4\vert\xi\vert^{2}-\log^{2}(1+\vert\xi\vert^{2\theta})} \leq \frac{\log^{2}(1+\vert\xi\vert^{2\theta})}{2\vert\xi\vert^{2}} \leq 1\]
for any $\xi \in \Omega_{h,2}$.
\hfill
$\Box$


In order to get Theorem 1.1 we need one more proposition in $\Omega_{h}$.
\begin{pro}\label{prop-3.4}
 Let $\theta>1/2$. Then, it holds that
 $$\int_{\Omega_h} |P_{1}e^{-a(\xi)t}\frac{\sin(t\vert\xi\vert)}{\vert\xi\vert}|^2d\xi \leq C\vert P_{1}\vert^{2}o(t^{-\frac{n}{2\theta}}), \quad(t \to \infty).$$
\end{pro}
{\it Proof.}\,Indeed,
\begin{align*}
\int_{\Omega_h} |e^{-a(\xi)t}\frac{\sin(t\vert\xi\vert)}{\vert\xi\vert}|^2d\xi
     &\leq \delta_{1}^{-2}\int_{\Omega_h} e^{-\log(1+\vert\xi\vert^{2\theta})t}d\xi = \delta_{1}^{-2}\int_{\Omega_h}(1+\vert\xi\vert^{2\theta})^{-t}d\xi \\
     & = \delta_{1}^{-2}\omega_{n}\int_{\delta_{1}}^{\infty}(1+r^{2\theta})^{-t}r^{n-1}dr \leq C t^{-1}\,(1+\delta_1^{2\theta})^{-t}
\end{align*}
for some constant~$C=C(n,\theta)>0$, where we applied Corollary~\ref{cor2.1} in the last inequality. This concludes the proof.
\hfill
$\Box$
\begin{rem}\,{\rm The decay rate stated in Proposition \ref{prop-3.4} can be drawn with a more precise fast decay rate, however, since the decay rate in Proposition \ref{proposition3.1} is essential, and the rate of decay in Proposition \ref{prop-3.4} can be absorbed into that of Proposition \ref{proposition3.1}, we just have employed such a style for simplicity.}
\end{rem}

Finally, Theorem \ref{main-theo} is a direct consequence of Propositions \ref{proposition3.1}, \ref{proposition3.2} and \ref{prop-3.4}. We shall draw its outline of proof in the case when $\theta \geq 3/4$.\\

{\it Outline of proof of Theorem 1.1.}\, Let $\theta>1/2$, and set
\[\nu_{\theta}(t,\xi) := P_{1}e^{-a(\xi)t}\frac{\sin(t\vert\xi\vert)}{\vert\xi\vert}.\]
Then, one can estimate as follows.
\[\int_{{\bf R}^{n}}\vert\hat{u}(t,\xi)-\nu_{\theta}(t,\xi)\vert^{2}d\xi\]
\begin{equation}\label{proof1}
= \left(\int_{\vert\xi\vert\leq\delta_{1}} + \int_{\Omega_{h}}\right)\vert\hat{u}(t,\xi)-\nu_{\theta}(t,\xi)\vert^{2}d\xi := I_{1}(t) + I_{2}(t).
\end{equation}
To begin with, by using Proposition \ref{proposition3.1} one has
\begin{equation}\label{proof2}
I_{1}(t) \leq C\Vert u_{1}\Vert_{1,1}^{2}t^{-\frac{n-(6-8\theta)_+}{2\theta}}\quad (t \gg 1),
\end{equation}
if~$n\ge2$ and~$\theta>1/2$, or~$n=1$ and~$\theta>5/8$. On the other hand,
\begin{equation}\label{proof2n1}
I_{1}(t) \leq C\Vert u_{1}\Vert_{1,1}^{2}t^{\frac1\theta-\frac32}\quad (t \gg 1),
\end{equation}
if~$n=1$ and~$\theta\in(1/2,5/8]$.
Secondary, from Propositions \ref{proposition3.2} and \ref{prop-3.4}
\begin{equation}\label{proof4}
I_{2}(t) \leq C\int_{\Omega_{h}}\vert\hat{u}(t,\xi)\vert^{2}d\xi + C\int_{\Omega_{h}}\vert\nu_{\theta}(t,\xi)\vert^{2}d\xi \leq C\,\|u_1\|^2\,e^{-\gamma t} + C\,|P_1|^{2}o(t^{-\frac{n}{2\theta}}), \quad(t \to \infty).
\end{equation}
The statement of Theorem 1.1 can be proved by combining \eqref{proof1}, \eqref{proof2} or \eqref{proof2n1}, and \eqref{proof4}.
\hfill
$\Box$

\section{Optimal asymptotic behavior}

In this section we study the optimality of various estimates of the integrals closely related with the leading terms obtained in previous sections.\\
We first prepare the following proposition in the large dimensional case.

\subsection{Optimal behavior  of the $L^2$-norm}
\begin{pro}\label{2.8}
Let  $n > 2$ and $\theta \geq 1/2$. Then there  exists $t_0 > 0$ such that  for $t \geq t_{0}$ it holds that
$$C^{-1} t^{-\frac{n-2}{2\theta}} \geq \int_{{\bf R}^{n}}{e^{-t \log(1+ |\xi|^{2\theta}) }\dfrac{|\sin(|\xi|t)|^2}{|\xi|^2}}d\xi \geq C t^{-\frac{n-2}{2\theta}}, $$
with  $C$ a positive constant depending only on   $n$ and $\theta$.
\end{pro}
{\it Proof.}\,  First, we may note that
\begin{align*}
M(t): & = \int_{{\bf R}^{n}} e^{-t\; \log(1+|\xi|^{2\theta})} \dfrac{|\sin(|\xi|t)|^2}{|\xi|^2} d\xi
\\
&=\omega_n \int_0^{\infty}e^{-t\; \log(1+r^{2\theta})}  r^{n-3}|\sin(rt)|^2   dr\\
& \geq\omega_n \int_0^{\infty}e^{-t\;r^{2\theta}}  r^{n-3}|\sin(rt)|^2   dr.\\
\end{align*}
Now we apply the change of variable $s=t^{1/2\theta}r$, for a fixed $t>0$, to arrive at
\begin{align*}
M(t) & \geq \omega_nt^{-\frac{n-2}{2\theta}}\int_0^{\infty}e^{-s^{2\theta}}s^{n-3}\sin^2\big( t^{\frac{2\theta-1}{2\theta}}s)ds.
\end{align*}
\noindent
For $\theta=1/2$ the result directly follows from this last estimate.  For $\theta >1/2$ we use the fundamental identity
\[2\sin^2 x = \left(1-\cos 2x\right),\]
to  obtain
\begin{align*}
M(t) &\geq \frac{1}{2}\omega_n t^{-\frac{n-2}{2\theta}} \int_0^{\infty}e^{-s^{2\theta}}s^{n-3}\Big(1-\cos(2t^{\frac{2\theta-1}{2\theta}}s)\Big)ds\\
&=\frac{1}{2} \omega_n t^{-\frac{n-2}{2\theta}} \big(A_{n,\theta}-F_{n,\theta}(t)\big),
\end{align*}
where
\[A_{n,\theta} = \displaystyle \int_0^{\infty}e^{-s^{2\theta}}s^{n-3}ds, \quad \displaystyle  F_{n,\theta}(t)=\int_0^{\infty}e^{-s^{2\theta}}s^{n-3}\cos\big(2t^{\frac{2\theta-1}{2\theta}}s\big)ds.\]
Due to  the fact $e^{-s^{2\theta}}s^{n-3}\; \in L^1({\bf R})$\; ($n>2$),  we can apply the Riemann-Lebesgue theorem to get
$$F_{n,\theta}(t) \to 0, \quad t \to \infty .  $$
\noindent
Then we conclude the existence of  $t_0 >0$ such that  $F_{n,\theta}(t) \leq \displaystyle{\frac{A_{n,\theta}}{2}}$ for all $t \geq t_0$. Thus, the half part of proposition is proved with $C=\dfrac{\omega_nA_{n,\theta}}{4}$.

Now we prove the estimate from above of the proposition.  Indeed,
\begin{align*}
M(t) &\leq \int_{{\bf R}^{n}} e^{-t\; \log(1+|\xi|^{2\theta})}\vert\xi\vert^{-2} d\xi \\
&=\omega_n \int_0^{\infty}e^{-t\; \log(1+r^{2\theta})}  r^{n-3} dr\\
&=\omega_{n}\int_{0}^{\infty}(1+r^{2\theta})^{-t}r^{n-3}dr \\
& \leq C_{n,\theta}\,t^{-\frac{n-2}{2\theta}}, \quad t \gg 1,
\end{align*}
where one has just used Corollary~\ref{cor2.1}. This estimate completes the proof of the proposition.
\hfill
$\Box$


\begin{pro}\label{proposition4.1} Let $n = 1$ and $\theta > 0$. Then it is true that
\[\int_{{\bf R}}(1+\vert\xi\vert^{2\theta})^{-t}\frac{\sin^{2}(t\vert\xi\vert)}{\vert\xi\vert^{2}}d\xi \sim t,\quad (t \gg 1).\]
\end{pro}

The proof of this  proposition is the same to the case $ \theta=1 $ which appears in a Lemma by Char\~ao-Ikehata \cite{Log-damping}. In  fact,  the expression $(1+r^{2\theta})^{-t}$, $ \theta  \ne 0 $
does not  change the proof of  the case  $ \theta=1 $. That is, the  result of the lemma is independent of
$ \theta$.

\vspace{0.3cm}

Next we also deal with the two dimensional case.

The following  proposition has  a version for the case $ \theta=1 $  in Char\~ao-Ikehata \cite{Log-damping} and its proof is also independent of $ \theta$.

\begin{pro}\label{proposition4.2}\,Let $n = 2$ and $\theta > 0$. Then it is true that
\[\int_{{\bf R}^{2}}(1+\vert\xi\vert^{2\theta})^{-t}\frac{\sin^{2}(t\vert\xi\vert)}{\vert\xi\vert^{2}}d\xi \sim \log t,\quad (t \gg 1).\]
\end{pro}

\vspace{0.2cm}
Finally, let us now prove Theorem \ref{main-theo2} at a stroke. \\

{\it Proof of Theorem \ref{main-theo2} completed.}\,It follows from the Plancherel theorem and triangle inequality, with some constant $C_{n} > 0$ one can get
\[C_{n}\Vert u(t,\cdot)\Vert \geq \vert P_{1}\vert\Vert (1+\vert\xi\vert^{2\theta})^{-\frac{t}{2}}\frac{\sin(t\vert\xi\vert)}{\vert\xi\vert}\Vert - \Vert\hat{u}(t,\cdot)-P_{1}(1+\vert\xi\vert^{2\theta})^{-\frac{t}{2}}\frac{\sin(t\vert\xi\vert)}{\vert\xi\vert}\Vert.\]
and
\[C_{n}\Vert u(t,\cdot)\Vert \leq \vert P_{1}\vert\Vert (1+\vert\xi\vert^{2\theta})^{-\frac{t}{2}}\frac{\sin(t\vert\xi\vert)}{\vert\xi\vert}\Vert + \Vert\hat{u}(t,\cdot)-P_{1}(1+\vert\xi\vert^{2\theta})^{-\frac{t}{2}}\frac{\sin(t\vert\xi\vert)}{\vert\xi\vert}\Vert.\]
These inequalities together with Theorem \ref{main-theo}, Propositions \ref{2.8}, \ref{proposition4.1} and \ref{proposition4.2} imply the desired estimates. This part is, nowadays, well-known (see \cite{I-14, IO-17}). We stress that~$6-8\theta<2$, for any~$\theta>1/2$, and this gives~$t^{-\frac{n-(6-8\theta)_+}{2\theta}} = \textit{o} (t^{-\frac{n-2}{2\theta}})$. In the case~$n=1$ and~$\theta\in(1/2,5/8)$, we also notice that $1/\theta-3/2<1/2$, for any~$\theta>1/2$.
\hfill
$\Box$


\subsection{Optimal behavior of  energy norm}

In this subsection, we get the optimal decay estimates of the total energy itself. Again the total energy $E_{u}(t)$ is defined by
\[2E_{u}(t) := \Vert u_{t}(t,\cdot)\Vert^{2} + \Vert\nabla u(t,\cdot)\Vert^{2}.\]
We will study the  topic based on the expression defined by \eqref{express1}.

First of all, the estimate of the part $\Vert\nabla u(t,\cdot)\Vert=\Vert\vert\xi\vert\hat{ u}(t, \cdot) \Vert$ is simple. It is sufficient to multiply the expression
(\ref{express1}) by $ | \xi |$ and to make estimates similarly to get decay rates of the $L^2$-norm of $u(t, \xi)$.
The result will be

\begin{equation}\label{decay-gradient}
C^{-1}t^{-\frac{n}{2\theta}}|P_1|^2 \leq \Vert \nabla u(t,\cdot)\Vert^2 \leq C I_1^{2}\; t^{-\frac{n}{2\theta}}\quad  (t \gg 1)
\end{equation}
with $C$ a positive constant depending on $n$ and $\theta$.

The delicate part is to get the precise estimate of the $L^2$-norm of $\hat{u}_t(t, \xi)$.  In order to do that, we take the time derivative of $\hat{u}(t, \xi)$ given by (\ref{express1}) to obtain  the following expression to the time derivative of the solution with $ r=| \xi |$,

\begin{align}\label{time-derivative}
\hat{u}_t(t,\xi) &= - P_{1}a(\xi)e^{-a(\xi)t}\frac{\sin(tr)}{r} + P_{1}e^{-a(\xi)t}\cos(tr) \nonumber\\
&- P_{1}\frac{\log^{2}(1+r^{2\theta})}{8r^{3}}\frac{a(\xi)}{\sqrt{(1-\eta_{2}g(r))^{3}}}e^{-a(\xi)t}\sin(tr)\nonumber\\
&+P_{1}\frac{\log^{2}(1+r^{2\theta})}{8r^{3}}\frac{1}{\sqrt{(1-\eta_{2}g(r))^{3}}}e^{-a(\xi)t}r\cos(tr)
\\
&- \left(\frac{A_{1}(\xi) - iB_{1}(\xi)}{b(\xi)}\right)a(\xi)e^{-a(\xi)t}\sin(b(\xi)t)
+ \Big(A_{1}(\xi) - iB_{1}(\xi)\Big)e^{-a(\xi)t}\cos(b(\xi)t)
\nonumber\\
&- tP_{1} a(\xi)e^{-a(\xi)t}\left(\frac{b(\xi)-r}{b(\xi)}\right)\cos(\mu(\xi)t) - tP_{1} e^{-a(\xi)t}\left(\frac{b(\xi)-r}{b(\xi)}\right)\mu(\xi)\sin(\mu(\xi)t) \nonumber \\
& P_{1} e^{-a(\xi)t}\left(\frac{b(\xi)-r}{b(\xi)}\right)\cos(\mu(\xi)t) \nonumber .
\end{align}
It should be remarked that \eqref{time-derivative} holds for small frequency parameters $\vert\xi\vert \ll 1$ and all $\theta>1/2$. A candidate to be a leading term as $t \to \infty$ of the velocity $u_t(t,\xi)$ is given by a simple form:
\begin{equation}\label{leading2}
P_{1}e^{-a(\xi)t}\cos(|\xi|t)
\end{equation}
where $a(\xi)=\displaystyle{\frac{\log(1+|\xi|^{2\theta})}{2}}$ and $P_1:= \displaystyle{\int_{{\bf R}^{n}}}u_{1}(x)dx$.

Our goal in this section is to get decay estimates in time to the remainder therms that appear in \eqref{time-derivative}, and are defined by the next $8$ functions which imply remainders with respect to the leading term \eqref{leading2}.
 \begin{itemize}
\item[$\bullet$]$ F_1(t, \xi)=- P_{1}\displaystyle{\frac{\log^{2}(1+r^{2\theta})}{8r^{3}}}\displaystyle{\frac{a(\xi)}{\sqrt{(1-\eta_{2}g(r))^{3}}}}e^{-a(\xi)t}\sin(tr)$,

\item[$\bullet$]$F_2(t, \xi)=P_{1}\displaystyle{\frac{\log^{2}(1+r^{2\theta})}{8r^{3}}}\displaystyle{\frac{1}{\sqrt{(1-\eta_{2}g(r))^{3}}}}e^{-a(\xi)t}r\cos(tr)$,

 \item[$\bullet$] $ F_3(t, \xi)=- \left(\displaystyle{\frac{A_{1}(\xi) - iB_{1}(\xi)}{b(\xi)}}\right)a(\xi)e^{-a(\xi)t}\sin(b(\xi)t)$,

 \item[$\bullet$] $F_4(t, \xi)= \Big(A_{1}(\xi) - iB_{1}(\xi)\Big)e^{-a(\xi)t}\cos(b(\xi)t)$,

 \item[$\bullet$] $F_5(t, \xi)= tP_{1} a(\xi)e^{-a(\xi)t}\left(\displaystyle{\frac{b(\xi)-r}{b(\xi)}}\right)\cos(\mu(\xi)t)$,

 \item[$\bullet$] $F_6(t, \xi)=tP_{1} e^{-a(\xi)t}\left(\displaystyle{\frac{b(\xi)-r}{b(\xi)}}\right)\mu(\xi)\sin(\mu(\xi)t)$,

 \item[$\bullet$] $F_7(t, \xi)= P_{1} e^{-a(\xi)t}\left(\displaystyle{\frac{b(\xi)-r}{b(\xi)}}\right)\cos(\mu(\xi)t)$,

 \item[$\bullet$] $F_8(t, \xi)= P_{1}a(\xi)e^{-a(\xi)t}\displaystyle{\frac{\sin(t|\xi|)}{|\xi|}}$,
 \end{itemize}
where $r:=|\xi|$.
 Then, the remainder term for  $\hat{u}_t(t,\xi)$ is given by
\begin{equation}\label{expressionder}
\hat{u}_t(t,\xi) -P_{1}e^{-a(\xi)t}\cos(tr) = \sum_{j = 1}^{8}F_{j}(t,\xi).
\end{equation}

We need the following lemmas in order to show that the remainder terms decay faster than the leading term (asymptotic profile).

\begin{lem}\label{lem-sen}
Let  $n \geq 1$ and $\theta \geq 1/2$. Then there  exists $t_0 > 0$ such that  for $t \geq t_{0}$ it holds that
$$Ct^{-\frac{n+4\theta-2}{2\theta}} \geq \int_{{\bf R}^{n}}{a(\xi)^2e^{-t \log(1+ |\xi|^{2\theta}) }\dfrac{|\sin(|\xi|t)|^2}{|\xi|^2}}d\xi \geq C^{-1} t^{-\frac{n+4\theta-2}{2\theta}}, $$
with  $C$ a positive constant depending only on   $n$ and $\theta$.
\end{lem}

{\it Proof:}  First place we note that
$$\log(1+r^{2\theta}) \leq r^{2\theta},$$
for all $r \geq 0$ and $\theta >1/2$.

Then, using the definition of $a(\xi)$ we can easily obtain the estimate from above.

\begin{align*}
\int_{{\bf R}^{n}}{a(\xi)^2e^{-t \log(1+ |\xi|^{2\theta}) }\dfrac{|\sin(|\xi|t)|^2}{|\xi|^2}}d\xi
&= \omega_n\int_0^{\infty} \frac{1}{4}\log^2(1+r^{2\theta})e^{-t \log(1+ r^{2\theta}) }\dfrac{|\sin(rt)|^2}{r^2}r^{n-1}dr\\
&\leq \omega_n\int_0^{\infty}\dfrac{r^{4\theta}}{4}(1+ r^{2\theta})^{-t} \dfrac{|\sin(rt)|^2}{r^2}r^{n-1}dr\\
&\leq \dfrac{\omega_n}{4}\int_0^{\infty}r^{n+4\theta-3}(1+ r^{2\theta})^{-t} dr\\
&\leq C_{n,\theta} \,t^{-\frac{n+4\theta-2}{2\theta}}, \quad t \gg 1,
\end{align*}
with a positive constant $C_{n,\theta}$, where we have just used Corollary \ref{cor2.1}.
\vspace{0.2cm}

Next we want to get the lower bound. We first observe that
$$\lim_{r \rightarrow 0}\frac{\log(1+r^{2\theta})}{r^{2\theta}}=1 .$$
Thus, there exists $\delta_0=\delta_0(\theta)>0$ such that
\begin{align}\label{lim-log}
1/2 \leq \frac{\log(1+r^{2\theta})}{r^{2\theta}} \leq 2
\end{align}
for all  $0< r \leq \delta_0 $. By using \eqref{lim-log} we may obtain  for
$t \geq \frac{1}{2\delta_{0}^{2\theta}}$
\begin{align*}
& I_{\sin}(t) :=
\int_{{\bf R}^{n}}{a(\xi)^2e^{-t \log(1+ |\xi|^{2\theta}) }\dfrac{|\sin(|\xi|t)|^2}{|\xi|^2}}d\xi\\
&=\omega_n\int_0^{\infty} \frac{1}{4}\log^2(1+r^{2\theta})e^{-t\log(1+ r^{2\theta})} \dfrac{|\sin(rt)|^2}{r^2}r^{n-1}dr \\
&\geq  \frac{\omega_n}{4}\int_0^{\delta_0} \frac{r^{4\theta + n -3}}{4}e^{-2t r^{2\theta}}\sin^2(rt)dr\\
&=  \frac{\omega_n}{16\theta}2^{-\frac{1}{2\theta}-\frac{4\theta+n-3}{2\theta}}t^{-\frac{4\theta + n-2}{2\theta}}
\int_0^{\sqrt{2t}\delta_0^\theta} s^{\frac{3\theta + n -2}{\theta}}e^{-s^2}\sin^2(2^{-\frac{1}{2\theta}}t^{\frac{2\theta-1}{2\theta}}s^{\frac{1}{\theta}} )ds.
\end{align*}

We use the identity $\sin^2x=\displaystyle{\frac{1}{2}}(1-\cos(2x))$ to get

\begin{align*}
& I_{\sin}(t) \geq \frac{\omega_n}{16\theta}2^{-\frac{1}{2\theta}-\frac{4\theta+n-3}{2\theta}}t^{-\frac{4\theta + n-2}{2\theta}}
\int_0^{\sqrt{2t}\delta_0^\theta} s^{\frac{3\theta + n -2}{\theta}}e^{-s^2} \frac{1}{2}ds\\
& - \frac{\omega_n}{16\theta}2^{-\frac{1}{2\theta}-\frac{4\theta+n-3}{2\theta}}t^{-\frac{4\theta + n-2}{2\theta}}
\int_0^{\sqrt{2t}\delta_0^\theta}  s^{\frac{3\theta + n -2}{\theta}}e^{-s^2}\frac{1}{2}\cos(2^{1-\frac{1}{2\theta}}t^{\frac{2\theta-1}{2\theta}}s^{\frac{1}{\theta}} )ds\\
& \geq  \frac{\omega_n}{32\theta}2^{-\frac{1}{2\theta}-\frac{4\theta+n-3}{2\theta}}t^{-\frac{4\theta + n-2}{2\theta}}
\int_0^{1}s^{\frac{3\theta + n -2}{\theta}}e^{-s^2} ds\\
& - \frac{\omega_n}{32\theta}2^{-\frac{1}{2\theta}-\frac{4\theta+n-3}{2\theta}}t^{-\frac{4\theta + n-2}{2\theta}}
\int_0^{\infty} s^{\frac{3\theta + n -2}{\theta}}e^{-s^2}\cos(2^{1-\frac{1}{2\theta}}t^{\frac{2\theta-1}{2\theta}}s^{\frac{1}{\theta}} )ds.\\
\end{align*}

We note that $s^{\frac{3\theta + n -2}{\theta}}e^{-s^2} \in L^1(0, \infty)$ because of $\theta>1/2$.

Then, we define the following  positive constant depending on $\theta$ and $n$
$$A_{n,\theta} = \int_0^{1}s^{\frac{3\theta + n -2}{\theta}}e^{-s^2} ds.$$

Moreover, by the Riemann-Lebesgue Lemma one has
$$ \lim_{t \rightarrow \infty}
\int_0^{\infty} s^{\frac{3\theta + n -2}{\theta}}e^{-s^2}\cos(2^{1-\frac{1}{2\theta}}t^{\frac{2\theta-1}{2\theta}}s^{\frac{1}{\theta}} )ds = 0$$
because of $\theta>1/2$. So, we can choose $t_0 \geq  \frac{1}{2\delta_{0}^{2\theta}}$ such that
$$\int_0^{\infty} s^{\frac{3\theta + n -2}{\theta}}e^{-s^2}\cos(2^{1-\frac{1}{2\theta}}t^{\frac{2\theta-1}{2\theta}}s^{\frac{1}{\theta}} )ds \leq  \frac{A_{n,\theta}}{2}$$
for all $t \geq t_0$.

Combining these results with the last estimate for $I_{sin}(t)$ we arrive at the desired estimate:
$$I_{\sin}(t) \geq \frac{\omega_n}{32\theta}2^{-\frac{1}{2\theta}-\frac{4\theta+n-3}{2\theta}}t^{-\frac{4\theta + n-2}{2\theta}}
(A_{n,\theta} - \frac{1}{2}A_{n,\theta}), \quad t \geq t_0, \; n\geq 1.$$
\hfill
$\Box$

\begin{lem}\label{em-cos}
Let  $n \geq 1$ and $\theta \geq 1/2$. Then there  exists $t_0 > 0$ such that  for $t \geq t_{0}$ it holds that
$$C^{-1} t^{-\frac{n}{2\theta}} \geq \int_{{\bf R}^{n}} {e^{-t \log(1+ |\xi|^{2\theta}) }|\cos(|\xi|t)|^2}d\xi \geq C t^{-\frac{n}{2\theta}}, $$
with a constant $C > 0$ depending only on $n$ and $\theta$.
\end{lem}

{\it Proof: } To get the upper estimate we have
\begin{align*}
I_{\cos}(t)& := \int_{{\bf R}^{n}}{e^{-t \log(1+ |\xi|^{2\theta}) }|\cos(|\xi|t)|^2}d\xi =
\int_{{\bf R}^{n}}(1+ |\xi|^{2\theta})^{-t} |\cos(|\xi|t)|^2d\xi\\
 &\leq  \omega_n \int_0^\infty(1+ |r|^{2\theta})^{-t} r^{n-1}dr \,\leq C\,t^{-\frac{n}{2\theta}}, \quad t \gg 1
\end{align*}
with a constant $C > 0$ depending only on $n$ and $\theta$ according to Corollary \ref{cor2.1}.

\vspace{0.2cm}
To get the estimate from below  we use \eqref{lim-log}. Then one  has
\begin{align*}
I_{\cos}(t)&=
 \omega_n\int_{0}^{\infty}  e^{-t(1+ |r|^{2\theta})}|\cos(rt)|^2 r^{n-1}dr  \geq  \omega_n\int_{0}^{1}  e^{-2tr^{2\theta}}\cos^2(rt) r^{n-1}dr\\
 & = \frac{ \omega_n}{\theta}2^{-\frac{n}{2\theta}}t^{-\frac{n}{2\theta}}\int_{0}^{\sqrt{2t}}  e^{-s^2} s^{\frac{n-\theta}{\theta}} \cos^2(2^{-\frac{1}{2\theta}}t^{\frac{2\theta - 1}{2\theta}}s^{\frac{1}{\theta}}) ds\\
 & \geq \frac{ \omega_n}{\theta}2^{-\frac{n}{2\theta}}t^{-\frac{n}{2\theta}}\int_{0}^{1}  e^{-s^2} s^{\frac{n-\theta}{\theta}} \cos^2(2^{-\frac{1}{2\theta}}t^{\frac{2\theta - 1}{2\theta}}s^{\frac{1}{\theta}}) ds,  \quad t \gg 1.
\end{align*}

Next, with the same argument as in Lemma \ref{lem-sen} via the Riemann-Lebesgue Lemma used to prove  the lower bound for $I_{\sin}(t)$, and the fact that $e^{-s^2}s^{\frac{n-\theta}{\theta}} \in L^1(0,1)$ one can obtain the estimate from below of this lemma.\\
\hfill $\Box$

\vspace{0.2cm}

Next we give various estimates to the functions $F_j(t, \xi)$, $j=1, \cdots, 7$, since the estimate for the term $F_8(t,\xi)$ is already given by Lemma \ref{lem-sen}.
\vspace{0.1cm}

To estimate $F_1(t,\xi)$ one notes that due to the limit in \eqref{flim}, if $\theta>3/4$, then there exists $\delta_0  \in (0,  1]$ such that

\begin{equation}\label{flim2}
\frac{\log^{2}(1+r^{2\theta})}{8r^{3}}  \leq 1/2
\end{equation}
for $0<r=|\xi| \leq \delta_0$.

While, to the function $g(r)$ defined in Section 3 such that
\[g(r) = \frac{\log^{2}(1+r^{2\theta})}{4r^{2}},\]
it happens that
$$\lim_{r \rightarrow 0} g(r)=0.$$
Therefore, there exists a number $\delta_1\in (0,1]$ such that
$$(1-\eta_2g(r))^3 \geq (1-g(r))^3 \geq (1/2)^3$$
for $0< r \leq \delta_1$.

Now, we define $\delta=min\{\delta_0\;,\;\delta_1\}$. Then, using the last  estimate and \eqref{flim2} in

$$\int_{|\xi| \leq \delta}|F_1(t,\xi)|^2 d\xi \leq |P_1|^2\int_{|\xi| \leq \delta} \Big|\frac{\log^{2}(1+r^{2\theta})}{8r^{3}}\frac{a(\xi)}{\sqrt{(1-\eta_{2}g(r))^{3}}}e^{-a(\xi)t}\sin(tr)\Big|^2d\xi,
$$
one can obtain that
$$\int_{|\xi| \leq \delta}|F_1(t,\xi)|^2 d\xi \leq C|P_1|^2\int_{|\xi| \leq 1} |\xi|^2e^{-2a(\xi)t}d\xi,
$$
with a constant $C>0$ because of the fact that $a(\xi)=\frac{1}{2}\log(1+|\xi|^{2\theta})\leq \frac{1}{2}|\xi|$ for
$|\xi| \leq 1$.

Using  Corollary \ref{cor2.1} it implies that
\begin{equation}\label{F1}
\int_{|\xi| \leq \delta}|F_1(t,\xi)|^2 d\xi \leq C |P_1|^2 t^{-\frac{n+2}{2\theta}}, \quad t \gg 1.
\end{equation}

If $\theta\in(1/2,3/4]$, then we use~\eqref{eq:thetaloss}, so that~$|P_1|^2 t^{-\frac{n+2}{2\theta}}$ is replaced by~$|P_1|^2 t^{-\frac{n+2-(6-8\theta)_+}{2\theta}}$ in the estimate for~$F_1$ (there is no need to distinguish~$n\geq2$ and~$n=1$ in this case). Again, we stress that~$6-8\theta<2$, for any~$\theta>1/2$. We estimate~$F_2$ as we did for~$F_1$.

The estimates for the other functions on low frequency zone are  similarly done by using the method to estimate functions $K_j(t,\xi), \; j=1,2,3$. In particular, when one estimates $F_6(t,\xi)$, it is necessary to use the inequality $|\mu(\xi)|^2=|\eta_1b(\xi)+(1-\eta_1)|\xi||^2\leq  10|\xi|^2$ for $|\xi| \leq 1$ and estimate similar to \eqref{delta}. The result is that there exists a number $\delta>0$ such that
\begin{equation}\label{Fj}
\int_{|\xi| \leq \delta}|F_j(t,\xi)|^2 d\xi \leq C |P_1|^2 t^{-\frac{n+2}{2\theta}}, \quad t \gg 1,
\end{equation}
for $j=4,5,6,7$.

Finally, the estimate for $F_3$ can be obtained by following the similar estimate to that of  $K_1$ in \eqref{K{1}}. The result is
\begin{equation}\label{F3}
\int_{|\xi| \leq \delta}|F_3(t,\xi)|^2 d\xi \leq C ||u_1||_{1,1}^2 t^{-\frac{n+2}{2\theta}}, \quad t \gg 1.
\end{equation}

Combining these estimates above one can get the following result on the difference between $\hat{u}_t(t,\xi)$ and the asymptotic profile given by \eqref{leading2}.

\begin{pro}\label{proposition4.4}\, Let $n \geq 1$ and $\theta >1/2$. Then, there exists a small constant $\delta \in (0,1]$ such that
\[\int_{\vert\xi\vert\leq \delta} \big|\hat{u}_t(t,\xi) -P_{1}e^{-a(\xi)t}\cos(tr)\big|^{2} d\xi
\leq C \Big[\;\big(\vert P_{1}\vert^{2}\,t^{(6-8\theta)_+} + \Vert u_{1}\Vert_{1,1}^{2}\big)t^{-\frac{n+2}{2\theta}} + |P_1|^2  t^{-\frac{n+4\theta-2}{2\theta}} \; \Big]\quad (t \gg 1),\]
with some generous constant $C=C_{n,\theta} > 0$ depending only on $\theta$ and $n$.
\end{pro}
Based on Proposition \ref{proposition4.4}, one can get the crucial result on the behavior for the time derivative of the solution.
\begin{pro}\label{proposition4.5}\, Let $n \geq 1$, $\theta >1/2$. Then, in the case of $1/2 < \theta \leq1$ it holds that
\[\int_{{\bf R}^n} \big|\hat{u}_t(t,\xi) -P_{1}e^{-a(\xi)t}\cos(tr)\big|^{2} d\xi
\leq C \Big[\;  \Vert u_{1}\Vert_{1,1}^{2} +  {\|u_1\|^2} \Big] t^{-\frac{n+4\theta-2}{2\theta }}\; \quad (t \gg 1),\]
and in the case of $\theta\geq1$ it is true that
\[\int_{{\bf R}^n} \big|\hat{u}_t(t,\xi) -P_{1}e^{-a(\xi)t}\cos(tr)\big|^{2} d\xi
\leq C \Big[\;  \Vert u_{1}\Vert_{1,1}^{2} +  {\|u_1\|^2} \Big] t^{-\frac{n+2}{2\theta }}\; \quad (t \gg 1),\]
with some generous constant $C=C_{n,\theta} > 0$ depending only on $\theta$ and $n$.
\end{pro}

{\it Proof:} According to the Proposition \ref{proposition4.4}, in order to prove the statement it suffices to get the estimates on the high frequency region $\Omega_h=\{\xi \in {\bf R}^{n}\,:\,\vert\xi\vert \geq \delta\}$. In fact,
\begin{align}\label{high-f}
\int_{\Omega_h} &\big|\hat{u}_t(t,\xi) -P_{1}e^{-a(\xi)t}\cos(tr)|^2 d\xi  \nonumber\\
&\leq
2\int_{ \Omega_h} \big|\hat{u}_t(t,\xi)|^2d\xi  +2|P_{1}|^2 \int_{|\xi| \geq \delta}e^{-2a(\xi)t}d\xi \nonumber\\
&\leq {\|u_1\|^{2}}\,e^{-\gamma t} +2|P_1|^2 \int_{\delta}^{\infty} ( 1+ r^{2\theta})^{-t}\omega_nr^{n-1}d\xi \\
&\leq {\|u_1\|^{2}}\,e^{-\gamma t} +2|P_1|^2 C \frac{(1+\delta^{2\theta})^{-t}}{t}\nonumber
\end{align}
which holds for $t \gg 1$, where we have just used Proposition \ref{proposition3.2} and Corollary~\ref{cor2.1}.
\hfill
$\Box$

{ \bf {  Proof of Theorem \ref{main-theo3}} :} The proof for lower bound of decay can be done by \eqref{decay-gradient} and the following estimate
$$||u_t(t,\cdot)|| \geq  || P_1e^{-a(\xi)t}\cos(tr)|| - ||u_t(t,\cdot) - P_1 e^{-a(\xi)t}\cos(tr)||
$$
combined with the estimates from below of Lemma \ref{em-cos} and Proposition \ref{proposition4.5}.

The estimate from above can be obtained using by Proposition \ref{proposition3.2},
the estimates   for the functions $F_j(t,\xi)$ $(j=1,\cdots,8)$ and the estimate from above of Lemma \ref{em-cos}.
\hfill
$\Box$

\par
\vspace{0.5cm}
\noindent{\em Acknowledgement.}
\smallskip
{\rm The work of the first author (R. C. CHAR\~AO) was partially supported by PRINT/CAPES - Process 88881.310536/2018-00 and the work of the third author (R. IKEHATA) was supported in part by Grant-in-Aid for Scientific Research (C)20K03682  of JSPS.}
}


{\rm

}


\end{document}